\documentclass[10pt,reqno]{amsart}

\usepackage{amssymb}
\usepackage{amsthm}
\usepackage{amsmath}
\usepackage{enumerate}

\usepackage[dvipdfm]{graphicx}
\usepackage{indentfirst}

\newtheorem{thm}{Theorem}
\newtheorem{cor}{Corollary}
\newtheorem{lem}{Lemma}
\theoremstyle{definition}

\newtheorem{remark}{Remark}
\renewcommand{\Re}{{\operatorname{Re}\,}}

\title{Subordination problems of Robertson functions}
\author[L.-M. Wang]{Li-Mei Wang}
\date{}

\begin{document}

\maketitle

\begin{abstract}
In the present paper, we are concerned with subordination problems related to $\lambda$-Robertson function. The radii of $\lambda$-spirallikeness and starlikeness of $\lambda$-Robertson function are also determined.

\end{abstract}

\footnote[0]{Primary 30C45, Secondary 30C62, 30C75}

\section{Introduction and Main Results}
Let $\mathbb{D}_{r}=\{z\in \mathbb{C}:\, |z|<r\}$ for $0<r\leq 1$ and $\mathbb{D}=\mathbb{D}_{1}$ be the unit disc. Let $\mathcal{A}$ be the family of functions $f$ analytic in $\mathbb{D}$, and $\mathcal{A}_{1}$ be the subset of $\mathcal{A}$ consisting of functions $f$ which are normalized by $f(0)=f'(0)-1=0$. 
A function $f\in\mathcal{A}$ is said to be subordinate to a function $F\in\mathcal{A}$ in $\mathbb{D}$ (in symbols $f\prec F$ or $f(z)\prec F(z)$) if there exists an analytic function $\omega(z)$ on $\mathbb{D}$ with $|\omega(z)|<1$ and $\omega(0)=0$, such that 
\[
f(z)=F(\omega(z))
\]
in $\mathbb{D}$. When $F$ is a univalent function, the condition $f\prec F$ is equivalent to $f(\mathbb{D})\subseteq F(\mathbb{D})$ and $f(0)=F(0)$.
 Let 
\[
\mathcal{P}_{\lambda}=\left\{p\in\mathcal{A} :\,  p(0)=1,\,\, \Re e^{-i \lambda}p(z)>0 \right\}.
\]
Here and hereafter we always suppose $-\pi/2<\lambda<\pi/2$. Note that $\mathcal{P}_{\lambda}$ is a convex and compact subset of $\mathcal{A}$ which is equipped with the topology of uniform convergence on compact subsets of $\mathbb{D}$.  Since $\mathcal{P}_{0}$ is the well-known Carath{\'e}odory class, we call $\mathcal{P}_{\lambda}$ \textit{the tilted Carath{\'e}odory class by angel $\lambda$}. Some characterizations and estimates of elements in $\mathcal{P}_{\lambda}$ are known (for a short survey, see \cite{Wang}).

For a function $f\in\mathcal{A}$, let 
\begin{equation}\label{Q(z)}
Q_{f}(z)=\frac{zf'(z)}{f(z)}
\end{equation}
and
\begin{equation*}
P_{f}(z)=1+\frac{zf''(z)}{f'(z)}.
\end{equation*}
It is worthwhile to note that
\[
Q_{f}(z)+\frac{zQ_{f}'(z)}{Q_{f}(z)}=P_{f}(z).
\]
These quantities are important for investigation of geometric properties of analytic functions. Next we will define two subclasses of analytic functions related to these two quantities.

A function $f\in\mathcal{A}_{1}$ is said to be a \textit{$\lambda$-spirallike function} (denoted by $f\in\mathcal{SP}_{\lambda}$) if
\[
Q_{f}\in \mathcal{P}_{\lambda}.
\]
Note that $\mathcal{SP}_{0}$ is precisely the set of starlike functions normally denoted by $\mathcal{S}^{*}$. Sprirallike functions were introduced and proved to be univalent by \v Spa\v cek \cite{Spacek} in 1932. 
For general references about spirallike functions, see e.g. \cite{Duren} or \cite{AhujaSilverman}.

A function $f\in\mathcal{A}_{1}$ is said to be a \textit{$\lambda$-Robertson function} if $zf'(z)\in\mathcal{SP}_{\lambda}$, i.e. 
\[
P_{f}\in \mathcal{P}_{\lambda}.
\] Let $\mathcal{R}_{\lambda}$ denote the set of these functions. Note that $\mathcal{R}_{0}$ is precisely the set of convex functions sometimes denoted by $\mathcal{K}$. Convex functions have been the subject of numerous investigations, among which the following result was proved by MacGregor \cite{MacGregor} in 1975.\\

\noindent{\bf Theorem A.}
\textit{\,\,Let $f\in \mathcal{K}$, then the subordination relation
\[
\frac{zf'(z)}{f(z)}\prec \frac{zf'_{0}(z)}{f_{0}(z)}
\]
holds, where $f_{0}(z)=z/(1-z)$}.\\

We are interested in more general subordination problems related to $\lambda$-Robertson functions. For this purpose, we first introduce some specific functions for convenience.

A distinguished member of $\mathcal{R}_{\lambda}$ is 
\begin{equation*}
f_{\lambda}(z)=\frac{(1-z)^{1-2e^{i\lambda}\cos\lambda}-1}{2e^{i\lambda}\cos\lambda -1}.
\end{equation*}
A simple calculation yields
\begin{equation}\label{q(lambda)}
Q_{\lambda}(z)=\frac{zf'_{\lambda}(z)}{f_{\lambda}(z)}=\frac{e^{2i\lambda}z}{1-z-(1-z)^{1+e^{2i\lambda}}},
\end{equation}

\begin{equation*}
P_{\lambda}(z)=1+\frac{zf''_{\lambda}(z)}{f'_{\lambda}(z)}=\frac{1+e^{2i\lambda}z}{1-z}
\end{equation*}
and
\begin{equation}\label{differential2}
Q_{\lambda}(z)+\frac{zQ_{\lambda}'(z)}{Q_{\lambda}(z)}=P_{\lambda}(z).
\end{equation}

In \cite{KimSrivastava}, Kim and Srivastava posed the open problem which is an extension of Theorem A whether
\[
\frac{zf'(z)}{f(z)}\prec \frac{zf_{\lambda}'(z)}{f_{\lambda}(z)}
\]
holds for $f\in\mathcal{R}_{\lambda}$ with general $\lambda$. In other words, if 
\begin{equation}\label{broit-bouquet}
q(z)+\frac{zq'(z)}{q(z)}\prec P_{\lambda}(z)
\end{equation}
in $\mathbb{D}$, then whether
\begin{equation*}
q\prec Q_{\lambda}
\end{equation*}
holds in $\mathbb{D}$ for general $\lambda$? Relation ($\ref{broit-bouquet}$) is a kind of Briot-Bouquet differential subordinations which have a surprising number of important applications in the theory of univalent functions. Many sources and references are given in \cite{MillerMocanu}.
\par In the present paper, we solve the above problem in a restricted disc and obtain the radii of spirallikeness and starlikeness for Robertson funtions as well.

\begin{thm}\label{theorem 1}
Let $q\in \mathcal{A}$ with $q(0)=1$ satisfy the differential subordination 
\begin{equation}\label{p(z)}
q(z)+\frac{zq'(z)}{q(z)}\prec P_{\lambda}(z)
\end{equation}
and $R_{1}(\lambda)$ be defined by
\begin{equation}\label{R_{1}}
R_{1}(\lambda)=\sup\{r<1\,: Q_{\lambda}(rz)\prec P_{\lambda}(rz)\,\, in \, \mathbb{D}\}.
\end{equation}
Then
\[
q(z)\prec Q_{\lambda}(z) 
\]
in $|z|<R_{1}(\lambda)$.
\end{thm}

By the discussion in Section 1, we can deduce the following corollary immediately from Theorem 1.
\begin{cor}\label{corollary 1}
Let $f\in \mathcal{R}_{\lambda}$, then
\[
\frac{zf'(z)}{f(z)}\prec \frac{zf_{\lambda}'(z)}{f_{\lambda}(z)}
\]
in $|z|<R_{1}(\lambda)$, where $R_{1}(\lambda)$ is given in \eqref{R_{1}}.
\end{cor}

\begin{remark}
The radius of $\lambda$-spirallikeness of $\lambda$-Robertson functions is at least $R_{1}(\lambda)$, since 
\[
\frac{zf'(z)}{f(z)}\prec \frac{zf'_{\lambda}(z)}{f_{\lambda}(z)}\prec P_{\lambda}(z) 
\]
in $\mathbb{D}_{R_{1}(\lambda)}$ for any $f\in\mathcal{R}_{\lambda}$.
\end{remark}

\begin{remark}

For $0\leq r<1$, let 
\begin{equation*}
\begin{split}
\psi_{\lambda}(r)
&=\max_{|z|=1}|P_{\lambda}^{-1}(Q_{\lambda}(rz))/r|\\
&=\max_{|z|=1}\left|\frac{1}{r} \frac{mrz-1+(1-rz)^m}{1-(1-rz)^m}\right|\\
\end{split}
\end{equation*}
where $m=1+e^{2i\lambda}$. $\psi_{\lambda}(r)$ is an increasing function defined on $[0,1)$ with $\psi_{\lambda}(0)=0$. By the definition of subordination, $R_{1}(\lambda)$ defined in ($\ref{R_{1}}$) could be expressed in terms of $\psi_{\lambda}$:
\[
R_{1}(\lambda)=\sup\{r<1\,: \psi_{\lambda}(r)<1\}=\psi_{\lambda}^{-1}(1).
\] 
\end{remark}

\begin{thm}\label{theorem 2}
Let $q\in \mathcal{A}$ with $q(0)=1$ satisfy the differential subordination
\begin{equation}
q(z)+\frac{zq'(z)}{q(z)}\prec P_{\lambda}(z),
\end{equation}
then
\[
q(R_{2}z)\prec P_{0}(z) 
\]
in $\mathbb{D}$, where
\begin{equation}\label{R_{2}}
R_{2}:=R_{2}(\lambda)=\frac{2}{\sqrt{4+2\sqrt{3}|\sin(2\lambda)|}}.
\end{equation}
\end{thm}

\begin{cor}\label{corollary 2}
The radius of starlikeness of $\lambda$-Robertson functions is at least $R_{2}(\lambda)$ given in \eqref{R_{2}}.
\end{cor}

Note that $R_{1}(0)=R_{2}(0)=1$, thus both Corollary 1 and Corollary 2 imply Theorem A. Note also that in \cite{AhujaSilverman}, Ahuja and Silverman posed the problem to find the radius of starlikeness for all $\lambda$-Robertson functions. Corollary \ref{corollary 2} implies that this radius is at least
\[
R_{2}=\min\{R_{2}(\lambda):\,-\pi/2<\lambda<\pi/2\}=R_{2}(\pi/4)=\sqrt{3}-1\approx 0.732.
\]
Libera and Ziegler in \cite{Libera-Ziegler} have shown that the radius of close-to-convexity for all $\lambda$-Robertson functions is approximately 0.99097524 and the radius of convexity is $\sqrt{2}/2$.

\section{Proofs of Resluts}

In order to obtain our main results, the following lemmas are required.

\begin{lem}[\cite{Wang}]\label{lemma1}
Let $p\in \mathcal{P}_{\lambda}$, then we have
\begin{equation*}
\left|p(z)-A(r)\right|\leq B(r),
\end{equation*}
where 
\[
A(r)=\frac{1+r^2e^{2i\lambda}}{1-r^2},\,\,
B(r)=\frac{2r\cos \lambda}{1-r^2}
\]
and $r=|z|<1$. Equality holds if and only if $p(z)=P_{\lambda}(xz)$ with $|x|=1$.
\end{lem}

\begin{lem}[\cite{MillerMocanu}, Lemma 2.2d]\label{lemma2}
Let $g(z)$ and $h(z)$ be in $\mathcal{A}$ with $g(0)=h(0)$. If $g\not\prec h$ in $\mathbb{D}$, then there exist two points $z_{0}$ with $|z_{0}|<1$ and $\eta_{0}$ with $|\eta_{0}|=1$ and $s\geq 1$ such that
\[
g(\mathbb{D}_{|z_{0}|})\subset h(\mathbb{D}),
\]
\[
g(z_{0})=h(\eta_{0})
\]
and
\[
z_{0}g'(z_{0})=s \eta_{0}h'(\eta_{0}).
\]
\end{lem}

The next lemma is due to Nunokawa \cite{Nunokawa}. We only quote the relevant part.
\begin{lem}[\cite{Nunokawa}]\label{lemma3}
Let $p(z)\in\mathcal{A}$ satisfy $p(0)=1$ and $p(z)\not=0$ in $\mathbb{D}$. If there exsits a point $z_{0}\in\mathbb{D}$ such that $\Re p(z)>0$ in $|z|<|z_{0}|$ and $p(z_{0})=i a$ where $a\in\mathbb{R}\setminus \{0\}$, then 
\[
\frac{z_{0}p'(z_{0})}{p(z_{0})}=ik
\]
where $k\geq (a+1/a)/2$ if $a>0$ and $k\leq -(a+1/a)/2$ if $a<0$. 
\end{lem}

\noindent\textit{Proof of Theorem 1.} For simplicity, we let $R=R_{1}(\lambda)$ and $p(z)=q(z)+zq'(z)/q(z)$, thus $p\in\mathcal{P}_{\lambda}$. If $q(z)\not\prec Q_{\lambda}(z)$ in $|z|<R$, then Lemma \ref{lemma2} implies the existences of $z_{0}$ with $|z_{0}|<R$, $\eta_{0}$ with $|\eta_{0}|=R$ and $s\geq 1$ such that
\begin{equation}{\label{M-M condition}}
\begin{split}
q(\mathbb{D}_{|z_{0}|})\subset Q_{\lambda}(\mathbb{D}),\\
q(z_{0})=Q_{\lambda}(\eta_{0}),\\
z_{0}q'(z_{0})=s \eta_{0}Q_{\lambda}'(\eta_{0}).\\
\end{split}
\end{equation}

Thus in view of  \eqref{differential2}, $(\ref{p(z)})$ and $(\ref{M-M condition})$, we have
\begin{equation}{\label{equality}}
\begin{split}
p(z_{0})
&=q(z_{0})+\frac{z_{0}q'(z_{0})}{q(z_{0})}\\
&=Q_{\lambda}(\eta_{0})+s\frac{\eta_{0}Q_{\lambda}'(\eta_{0})}{Q_{\lambda}(\eta_{0})}\\
&=sP_{\lambda}(\eta_{0})+(1-s)Q_{\lambda}(\eta_{0}).
\end{split}
\end{equation}

Therefore ($\ref{equality}$) and Lemma $\ref{lemma1}$ show that
\begin{eqnarray*}
|p(z_{0})-A(R)|
&=&
|s(P_{\lambda}(\eta_{0})-A(R))-(1-s)(Q_{\lambda}(\eta_{0})-A(R))|\\
&\geq &
|s(P_{\lambda}(\eta_{0})-A(R))|-(s-1)|(Q_{\lambda}(\eta_{0})-A(R))|\\
&\geq&
sB(R)-(s-1)B(R)\\
&=&
B(R)
\end{eqnarray*}
which contradicts with $p(z)\in\mathcal{P}_{\lambda}$. Therefore we get the assertion.
\qed

\noindent\textit{Proof of Theorem 2.} For simplicity, we let $R=R_{2}(\lambda)$ and $p(z)=q(z)+zq'(z)/q(z)$, thus $p\in\mathcal{P}_{\lambda}$. If $q(Rz)\not\prec P_{0}(z)$ in $\mathbb{D}$, it follows from $q(0)=1$ that there esists a point $z_{0}\in\mathbb{D}$ such that $\Re q(Rz)>0$ for $|z|<|z_{0}|$ and $q(Rz_{0})=ia$ where $a\in\mathbb{R}\setminus \{0\}$, then by Lemma \ref{lemma3}, we have
\[
\frac{Rz_{0}q'(Rz_{0})}{q(Rz_{0})}=ik,
\]
where $k\geq (a+1/a)/2$ if $a>0$ and $k\leq -(a+1/a)/2$ if $a<0$. Therefore
\[
p(Rz_{0})=q(Rz_{0})+\frac{Rz_{0}q'(Rz_{0})}{q(Rz_{0})}=ia+ik,
\]
which implies $p(Rz_{0})\in\Omega$ since $|a+k|\geq \sqrt{3}$, where $\Omega=\{it,|t|\geq\sqrt{3}\}$. Next we will show that 
\[
p(\mathbb{D}_{R})\cap \Omega=\varnothing,
\]
which contradicts the above assertion. Since $p\in\mathcal{P}_{\lambda}$, it is sufficient to prove for functions $P_{\lambda}(z)$. Suppose that there is a point $z_{1}\in\mathbb{D}$ such that $P_{\lambda}(z_{1})=it_{0}$ with $|t_{0}|\geq\sqrt{3}$, then a simple calculation gives that
\[
z_{1}=\frac{it_{0}-1}{it_{0}+e^{2i\lambda}}.
\]
Hence 
\[
|z_{1}|^2=\frac{t_{0}^2+1}{t_{0}^2+1-2t_{0}\sin(2\lambda)}\geq\frac{4}{4+2\sqrt{3}|\sin(2\lambda)|}=R
\]
since $|t_{0}|\geq\sqrt{3}$. Therefore $P_{\lambda}(\mathbb{D}_{R})\cap\Omega=\varnothing$. The proof is completed.
 \qed

\[
\text{Acknowledgements}
\]
The author is grateful to Professor Toshiyuki Sugawa for his guidance and useful discussions during the preparation of this paper.

Division of Mathematics\\
Graduate School of Information Sciences\\
Tohoku University, Sendai\\
980-8579 JAPAN \\
e-mail: rime@ims.is.tohoku.ac.jp\\
\quad \quad \quad wangmabel@163.com

\end{document}